\documentclass[12pt,twoside]{article}
\usepackage[ansinew]{inputenc}
\usepackage{amsfonts}
\usepackage{latexsym}
\usepackage{amsmath}
\usepackage{amssymb}
\usepackage{eepic}
\usepackage{graphicx}
\usepackage{float}
\usepackage{pstricks}
\usepackage{epic,eepic}

\usepackage{xcolor}
\usepackage{tikz}
\usetikzlibrary{arrows,shapes,chains}

\newcommand{\R}{\mathbb{R}}

\newcommand{\C}{\mathbb{C}}
\newcommand{\N}{\mathbb{N}}
\newcommand{\Z}{\mathbb{Z}}

\newtheorem{defin}{Definition}[section]

\newtheorem{theorem}[defin]{Theorem}

\newtheorem{exa}{Example}
\newenvironment{example}{\begin{exa}\rm}{\end{exa}}

\newenvironment{proof}
{\noindent{\it Proof.}}{\hfill $\Box$\par\vspace{2.5mm}}

\newtheorem{que}{Question}

\newtheorem{pro}{Problem}

\numberwithin{equation}{section}

\makeatletter
\renewcommand{\ps@myheadings}{%
\renewcommand{\@evenhead}%
{{\rm\thepage}\hfil{\sc  J.~Heittokangas and Z.-T.~Wen}\hfil}%
\renewcommand{\@oddhead}%
{\hfil{{\sc Zeros of exponential sums in critical strips}\hfil{\rm\thepage}}}%
\renewcommand{\@evenfoot}{}%
\renewcommand{\@oddfoot}{}%
}\makeatother 

\title{\bf\Large The asymptotic number of zeros of exponential sums in critical strips}
\author{Janne Heittokangas
and Zhi-Tao Wen
   }
\date{}
\begin{document}
\maketitle

\begin{abstract}
Normalized exponential sums are entire functions of the form
    $$
    f(z)=1+H_1e^{w_1z}+\cdots+H_ne^{w_nz},
    $$
where $H_1,\ldots, H_n\in\C$ and $0<w_1<\ldots<w_n$.
It is known that the zeros of such functions are in finitely many
vertical strips $S$. The asymptotic number of the zeros in the union of all these strips was found by R.~E.~Langer already in 1931.
In 1973, C.~J.~Moreno proved that there are zeros arbitrarily
close to any vertical line in any strip $S$, provided that
$1,w_1,\ldots,w_n$ are linearly independent over the
rational numbers. In this study the asymptotic number of zeros in each
individual vertical strip is found by relying on R.~J.~Backlund's lemma,
which was originally used to study the zeros of the Riemann $\zeta$-function. As a counterpart
to Moreno's result, it is shown that almost every vertical line
meets at most finitely many small discs around the zeros of $f$.


\end{abstract}

\section{Introduction}
An exponential sum is an entire function of the form
    \begin{equation}\label{exp.eq}
    g(z)=F_1e^{\lambda_1z}+\cdots+F_me^{\lambda_mz},
    \end{equation}
where $\lambda_1, \ldots, \lambda_m\in\C$ are distinct constants called \emph{frequencies} or \emph{leading coefficients} of $g$,
and the coefficients $F_1,\ldots,F_m$ are polynomials.
The trigonometric functions $\sin z= \frac{1}{2i}\left(e^{iz}-e^{-iz}\right)$ and $\cos z= \frac{1}{2}\left(e^{iz}+e^{-iz}\right)$ are typical
examples of such functions. The exponential terms $e^{\lambda_1z},\ldots,e^{\lambda_mz}$ form a fundamental solution base to a linear differential equation with constant coefficients determined by the characteristic equation
	$$
	(r-\lambda_1)\cdots (r-\lambda_m)=0.
	$$
Hence, if $F_1,\ldots,F_m\in\C$, then the function $g$ in \eqref{exp.eq} is the generic solution of this differential equation.

We consider the zero distribution of exponential sums $g(z)$ in \eqref{exp.eq} in the case when all leading coefficients $\lambda_j$ are in one line, and $F_j$'s are complex constants.
By appealing to a rotation, if necessary, we may suppose that $\lambda_j$'s are all real and ordered such that $\lambda_1<\lambda_2<\ldots<\lambda_m$.
From the zero distribution point of view, we may pull the zero-free term $F_1e^{\lambda_1z}$ out as a common factor, and consider the zeros of
its multiplicative exponential sum. In other words, we may suppose that the
exponential sum under consideration is in the normalized form
    \begin{equation}\label{normal.eq}
    f(z)=1+H_1e^{w_1z}+\cdots+H_ne^{w_nz}
    \end{equation}
to begin with, where $0<w_1<\ldots<w_m$, $H_j$'s are non-zero complex constants, and $n\leq m-1$. Here we may suppose that $n\geq 2$, since in the case $n=1$ it is easy to see that the zeros of $f$ are all on one vertical line.

In 1917, Tamarkin found in his PhD thesis, published in Russian, that the zeros of exponential sums of the form \eqref{normal.eq} are confined to a vertical strip in the complex plane. An English translation of Tamarkin's
thesis appeared in 1928, see \cite{Tamarkin}. Wilder \cite{Wilder} gave an estimate for the number of
zeros in a rectangle with vertical sides on the boundary of a strip found by Tamarkin. An improvement of Wilder's result by Langer \cite{Langer} shows that the number $n(r)$ of zeros of a function $f$ of the form \eqref{normal.eq} in such a rectangle of height $r$ is subject to the bound
    \begin{equation}\label{nbdd.eq}
    \bigg|n(r)-\frac{w_n}{2\pi}r\bigg|\leq n.
    \end{equation}

In the case of polynomial coefficients, the zeros of exponential sums
are no longer in rectangular strips. A simple example is $g(z)=e^z-z$, whose zeros
are precisely on the curves $y=\pm \sqrt{e^{2x}-x^2}\sim \pm e^x$
\cite[Example~2.1]{HITW}. If $z=re^{i\theta}=r\cos\theta+ir\sin\theta=x+iy$ is on such curves, then $|\theta|\sim \frac{\pi}{2}-\frac{\log r}{r}$.
In general, the results by P\'{o}lya \cite{Poly, Poly2} and Schwengeler \cite{Sch} show that there exist finitely many angles $\theta^*$ such
that all but finitely many of the zeros of exponential sums of the form \eqref{normal.eq}, with polynomial coefficients, are located in the logarithmic strips of the form
	\begin{equation*}
	\Lambda(\theta^*,c)
	=\left\{re^{i\theta}: r>1,\, |\theta-\theta^*|
	<c\frac{\log r}{r}\right\},
	\end{equation*}
where $c>0$ is large enough. In addition, the number $n(r)$ of the zeros
of $f$ in the disc $|z|<r$ is
    \begin{equation}\label{PS}
    n(r)=\frac{w_m}{\pi}r+O(1).
    \end{equation}
See \cite{HITW, Heittokangas-Wen} for a more thorough historical review
and for improvements.

By combining the results of Moreno \cite{Moreno} and More et al.~\cite{Mora-Sepulcre-Vidal}, we obtain the following statement:
If the frequencies $w_1, \ldots, w_n$ of the normalized $f$ in \eqref{normal.eq} are linearly independent over the rational numbers,
then the zeros of $f$ are in vertical strips $S$, and the projections of the zeros of $f$ on $\R$ form dense subsets of the intervals induced by the
intersections of the strips $S$ with $\R$. Following Moreno \cite{Moreno}, we may express this by saying that there are zeros near
every line in every strip $S$. Alternatively,  for every vertical line $\textnormal{Re}\,(z)=\sigma$ in $S$ and for every $\varepsilon>0$
a point $\zeta=a+ib$ can be found such that $\sigma-\varepsilon<a<\sigma+\varepsilon$ and $f(\zeta)=0$. A precise statement is as follows.

\bigskip
\noindent
\textbf{Theorem A}
\emph{
Let $f(z)$ be an exponential sum of the normalized form \eqref{normal.eq}, whose leading coefficients $0<w_1<\ldots<w_n$ are linearly independent over rational numbers. Then an open interval $(\sigma_0, \sigma_1)$ is
contained in
	$$
    R_f:=\overline{\{\textnormal{Re}\,(z): f(z)=0\}}
    $$
if and only if the $n+1$ inequalities
    \begin{equation*}
    1\leq\sum_{j=1}^m|H_j|e^{w_j\sigma},
    \quad |H_k|e^{w_k\sigma}\leq 1+\sum_{j=1, j\neq k}^m|H_j|e^{w_j\sigma},\quad k=1,\ldots,n,
    \end{equation*}
are satisfied for any $\sigma\in(\sigma_0, \sigma_1)$.
}
\bigskip

As Moreno observed, the zeros of exponential sums might be in several vertical strips. In Langer's result \eqref{nbdd.eq} the zeros are considered to
be in one wide strip that in fact might include several narrower strips.
So far, none of the results in the existing literature seems say anyhing about the precise number of zeros in each individual strip. The main result in this paper, see Theorem~\ref{Main.theorem} below, aims to fill in this gap. In the next section, in addition to stating Theorem~\ref{Main.theorem}, we will also shed some more light into the zero distribution of exponential sums in vertical strips.


\section{Zero distribution of exponential sums}

Let $f(z)$ be an exponential sum of the form \eqref{normal.eq}.
If $z=x+iy$ satisfies
    \begin{equation}\label{mineq.eq}
    |H_k|e^{w_kx}> \sum_{j\neq k}^n|H_j|e^{w_jx},\quad k=0,1,\ldots,n,
    \end{equation}
where $H_0=1$ and $w_0=0$, then $f$ has no zeros in region
    $$
    G:=\{ z=x+iy: x~\text{satisfies}~\eqref{mineq.eq},\, -\infty<y<\infty\}.
    $$
This follows from the inequalities
    \begin{equation}\label{zerofree.eq}
    |f(z)|\geq |H_k|e^{w_kx}-\left( 1+\sum_{j=1, j\neq k}^m|H_j|e^{w_jx}\right)>0.
    \end{equation}
Any such vertical strip $G$ is called a \emph{zero-free region} of $f$.
The zero-free region of $f$ to the extreme left is determined by the
term~$1$, and the one to the extreme right
by $H_ne^{w_nz}$. They both have unbounded width, while the rest of the zero-free regions, if any, have bounded width.
Since $f$ has $n+1$ terms, it is obvious that $f$ has at most $n+1$
zero-free regions.

\begin{example}
Let $f(z)=1+e^z+e^{2z}$. Then the two sets
	$$
	\left\{z: \textnormal{Re}\,(z)<\log\frac{\sqrt{5}-1}{2}\right\}
	\quad\textnormal{and}\quad
	\left\{z: \textnormal{Re}\,(z)>\log\frac{\sqrt{5}+1}{2}\right\}
	$$
are the only zero-free regions of $f$ although $f$ has three terms.
\end{example}

In each zero-free region, precisely one inequality
in \eqref{mineq.eq} holds. The exponential term $H_je^{w_jz}$, whose
modulus is strictly greater than the sum of others,
is called the \emph{dominating term} in this region.
The boundary line of each zero-free domain of $f$ is denoted by $L_k$ and is determined by the equality
    \begin{equation}\label{boundary.eq}
    |H_k|e^{w_kx}= \sum_{j\neq k}^m|H_j|e^{w_jx}.
    \end{equation}
It is easy to see that there exist at most $2n$ boundary lines. The closed area between two consecutive zero free regions
with dominating terms $H_je^{w_jz}$ and $H_ke^{w_kz}$ for $j\neq k$
is called a \emph{critical strip} of $f$, and is denoted by $\Lambda(j,k)$.
Every critical strip is a closed area and contains
its vertical boundary lines. All zeros of $f$ lie in the critical strips,
and the zeros may lie on the vertical boundary lines of a critical strip, as is seen in the next example.

\begin{example}
The zeros of the exponential sum
	$$
	f(z)=6-5e^z+e^{2z}=(e^z-2)(e^z-3)
	$$
lie on two lines $\textnormal{Re}\,(z)=\log 2$ and $\textnormal{Re}\,(z)=\log 3$.
The zero free regions of $f$ are $\{z: \textnormal{Re}\,(z)>\log 6\}$, $\{z: \log 2<\textnormal{Re}\,(z)<\log 3\}$ and
$\{z: \textnormal{Re}\,(z)<0\}$. Thus the critical strips of $f$ are $\Lambda(1,2):=\{z:\log 3\leq \textnormal{Re}\,(z)\leq\log 6\}$
and $\Lambda(0,1):=\{z: 0\leq \textnormal{Re}\,(z)\leq\log 2\}$.
The lines $\textnormal{Re}\,(z)=\log 2$ and $\textnormal{Re}\,(z)=\log 3$ are  boundary lines of $\Lambda(0,1)$
and $\Lambda(1,2)$, respectively, see Figure~\ref{domain.png}.
    \begin{figure}[H]\label{zero}
    \begin{center}
    \includegraphics[scale=0.2]{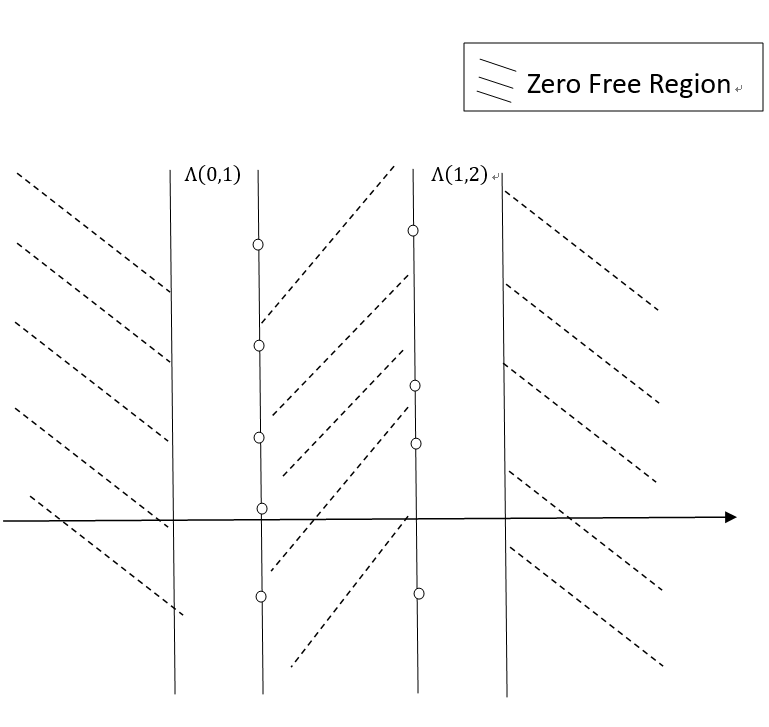}
    \end{center}
    \caption{Zero-free regions and critical strips.}\label{domain.png}
    \end{figure}
\end{example}

The following main result reveals the asymptotic number of zeros of exponential sums in each individual critical strip.

\begin{theorem}\label{Main.theorem}
Let $f$ be an exponential sum of the form \eqref{normal.eq}, where the coefficients are constants and
the leading coefficients satisfy $0<w_1<\ldots<w_n$. Then all zeros of $f$ are in finitely many critical strips $\Lambda(j,k)$.
Moreover, let $R$ be any rectangle cut from a critical strip $\Lambda(j,k)$
by two horizontal lines $\textnormal{Im}\,(z)=y_1$ and $\textnormal{Im}\,(z)\, =y_2$ with $y_2-y_1=r>0$. Then the number $n(r,\Lambda_{jk})$
of zeros of $f$ in $R$ satisfies
    \begin{equation}\label{n.eq}
    n(r,\Lambda_{jk})=\frac{|w_j-w_k|}{2\pi}r+O(1).
    \end{equation}
\end{theorem}

Let $R_1$ (resp.~$R_2$) be a rectangle cut from a critical strip $\Lambda(j,k)$
by the horizontal lines $\textnormal{Im}\,(z)=0$ and $\textnormal{Im}\,(z)\, =r>0$ (resp.~$\textnormal{Im}\,(z)=0$ and $\textnormal{Im}\,(z)\, =r+1$).
As an immediate consequence of the proof of \eqref{n.eq}, the number of zeros in $R_2\setminus R_1$ is asymptotically equal to $\frac{|w_j-w_k|}{2\pi}$ as $r\to\infty$.
For an analogous result in the lower half-plane, the horizontal lines $\textnormal{Im}\,(z)\, =r$ and $\textnormal{Im}\,(z)\, =r+1$ need to be replaced with $\textnormal{Im}\,(z)\, =-r$ and $\textnormal{Im}\,(z)\, =-r-1$, respectively.
This gives us more information about the distribution of zeros of $f$ in the critical strips.

Theorem A leads us to asking whether some
vertical lines in the zero-rich strips of $f$ could avoid the zeros of $f$? In the next example
we will see that this is indeed the case in the sense that almost every vertical line meets at
most finitely many small discs around the zeros of $f$.
This claim does not violate Theorem~A as the radius of these discs tends to zero.

\begin{example}
Let $f$ be as in Theorem~\ref{Main.theorem}, and let $\{z_n\}$ denote
the sequence of zeros of $f$ listed according to multiplicities and
ordered with respect to increasing modulus.
Let $\{D_n\}$ be the collection of Euclidean discs
    $$
    D_n:=\left\{z: |z-z_n|<r_n\right\},
    \quad r_n:=(1+|z_n|)^{-1}\log^{-2}(e+|z_n|).
    $$
We will show that the set $C\subset\R$ of values
$c$ for which the vertical line $\textnormal{Re}\,(z)=c$ meets infinitely
many discs $D_n$ has linear measure zero.

Let $n(r)$ denote the number of points $z_n$ in the disc $|z|<r$. By
Riemann-Stieltjes integration and \eqref{PS}, we obtain
	$$
	\sum_{n=1}^\infty r_n=\int_0^\infty \frac{dn(t)}{(1+t)\log^2(e+t)}
	\leq 3\int_0^\infty \frac{n(t)\, dt}{(1+t)^2\log^2(e+t)}<\infty.
	$$
Thus, for any $\varepsilon>0$ there exists
an $N(\varepsilon)\in\N$ such that $\sum_{n\geq N(\varepsilon)}2r_n<\varepsilon$.
Let $C_\varepsilon\subset\R$ denote the set of values $c$ for which the line $\textnormal{Re}\,(z)=c$
meets at least one of the discs $D_n$ with $n\geq N(\varepsilon)$. Then $C_\varepsilon$
has linear measure $<\varepsilon$. Since $C$ is contained in all sets $C_\varepsilon$,
$\varepsilon>0$, it follows that $C$ has linear measure zero.
\end{example}

The proof of Theorem~\ref{Main.theorem} in Section~\ref{proof-sec} is based on Backlund's lemma,
which will be discussed in Section~\ref{Backlund-sec}.


\section{Backlund's lemma}\label{Backlund-sec}

A slightly weaker version of the lemma below is proved in the PhD
thesis of Schwengeler \cite[pp.~34--39]{Sch}. Some years earlier, Backlund had used similar ideas in his considerations regarding the zeros of the Riemann $\zeta$-function \cite{Backlund}. Although at first glance
the lemma seems to be independent of this paper,
it is in a key role in proving Theorem~\ref{Main.theorem}.
The original proof in \cite{Sch} is hard to read due
to lack of details, lack of assumptions, and confusing
notation. Hence a complete proof is given for the convenience of the reader.

\bigskip
\noindent
\textbf{Backlund's lemma.}
\emph{Suppose that $F$ is an entire function, and that $z_1, z_2\in\C$ are distinct points such that $F$ has no zeros on the closed interval
$[z_1,z_2]$. Then
	$$
	\left|\textnormal{Re} \left\{\frac{1}{2\pi i}\int^{z_2}_{z_1}\frac{F'(z)}{F(z)}\,dz\right\}\right|
	\leq \frac{\max_{|\zeta|\leq R}\log|F(z_1+\zeta)|-\log |F(z_1)|}{2\log (R/T)}+\frac{1}{2},
	$$
where $R>T=|z_2-z_1|$.}

\bigskip
\begin{proof}
By appealing to a suitable change of variable, we may suppose that
$z_1,z_2$ have the same real part.
Without loss of generality, we may suppose that $\textnormal{Im}(z_2)>\textnormal{Im}\,(z_1)$. Then
$T=|z_2-z_1|=\textnormal{Im}\,(z_2)-\textnormal{Im}\,(z_1)>0$.
Denote the points between $z_1$ and $z_2$ by $z'=z_1+i\eta$, where $0\leq \eta\leq T$.

Since $F$ is entire and has no zeros on the interval $[z_1,z_2]$, we may use the fundamental
theorem of calculus to conclude that
	$$
	\int_{z_1}^{z_2}\frac{F'(z)}{F(z)}\, dz=\int_{z_1}^{z_2}\,d\log F(z)=\log F(z)\bigg|^{z_2}_{z_1},
	$$
from which
	\begin{equation}\label{RI.eq}
	\text{Re} \left\{\frac{1}{2\pi i}\int^{z_2}_{z_1}\frac{F'(z)}{F(z)}\,dz\right\}=\frac{1}{2\pi}
	\text{Im} \left\{\log F(z)\bigg|^{z_2}_{z_1}\right\}
	\end{equation}
holds for any fixed branch of the logarithm.
Set $w=F(z)=\rho e^{i\varphi}$, so that $\log F(z)=\log\rho+i\varphi$.
The argument of $F(z_j)$ is denoted by $\varphi_{z_j}$ for $j=1,2$. Then \eqref{RI.eq} yields
	\begin{equation}\label{Im.eq}
	\text{Re} \left\{\frac{1}{2\pi i}\int^{z_2}_{z_1}\frac{F'(z)}{F(z)}\,dz\right\}=
    \frac{1}{2\pi}(\varphi_{z_2}-\varphi_{z_1}).
	\end{equation}
If $z'=z_1+i\eta$ moves from $z_1$ to $z_2$, then its image
$F(z')$ moves along a curve $C$ from $F(z_1)$ to $F(z_2)$ as in Figure~\ref{FBC.pdf}, where $\Gamma$
is a line through the origin and orthogonal to the vector from the origin to the point $F(z_1)\neq 0$.
Let $s$ be the number of intersection points of the curve $C$ with the line $\Gamma$.
It is easy to see that
	\begin{equation}\label{number.eq}
	|\varphi_{z_2}-\varphi_{z_1}|\leq \pi s+\pi.
	\end{equation}

Using the Taylor series expansion $F(z)=\sum_{n=0}^\infty \frac{F^{(n)}(z_1)}{n!}(z-z_1)^n$, we obtain
	$$
	F(z_1+i\eta)=F(z_1)+C_1\eta+C_2\eta^2+\cdots+C_n\eta^n+\cdots.
	$$	
Following Schwengeler \cite{Sch}, we denote $\overline{F}(z)=\overline{F(\overline{z})}$. Then $\overline{F}$
is an entire function, and $\overline{F}(\overline{z})=\overline{F(z)}$, so that
	$$
    \overline{F}(\overline{z}_1-i\eta)=\overline{F}(\overline{z}_1)+\overline{C}_1\eta+\overline{C}_2\eta^2+\cdots+\overline{C}_n\eta^n+\cdots.
	$$	
Since $F(z_1)\neq 0$, we may set $a_j=C_j/F(z_1)$. Then $\overline{a}_j=\overline{C}_j/\overline{F(z_1)}
=\overline{C}_j/\overline{F}(\overline{z}_1)$, and we have
	\begin{eqnarray*}
	\frac{F(z_1+i\eta)}{F(z_1)} &=& 1+a_1\eta+a_2\eta^2+\cdots+a_n\eta^n+\cdots,\\
    \frac{\overline{F}(\overline{z}_1-i\eta)}{\overline{F}(\overline{z}_1)}
    &=&1+\overline{a}_1\eta+\overline{a}_2\eta^2+\cdots+\overline{a}_n\eta^n+\cdots.
    \end{eqnarray*}
Define $\Phi:\R\to\R$ by
	\begin{equation}\label{Phi.eq}
    \begin{split}
	\Phi(\eta)&=\frac{1}{2}\left(\frac{F(z_1+i\eta)}{F(z_1)}+
	\frac{\overline{F}(\overline{z}_1-i \eta)}{\overline{F}(\overline{z}_1)}\right)\\
    &=1+\frac{a_1+\overline{a}_1}{2}\eta+\frac{a_2+\overline{a}_2}{2}\eta^2
    +\frac{a_3+\overline{a}_3}{2}\eta^3+\cdots.
    \end{split}
	\end{equation}
We may analytically continue $\Phi$ into an entire function $\Phi(z)$ simply by replacing $\eta$ with $z$.
This follows from the representation
    $$
    \Phi(z)=\frac{1}{2}\left(\frac{G(z)}{F(z_1)}
    +\frac{\overline{G}(z)}{\overline{F}(\overline{z}_1)}\right),
    $$
where $G(z)=F(z_1+iz)$ is an entire function. Alternatively, we may see that $\Phi(z)$ is entire by
noticing that the power series in \eqref{Phi.eq}, with $\eta$ replaced by $z$, has infinite radius
of convergence.

    \begin{figure}[H]
    \begin{center}
    \includegraphics[scale=0.4]{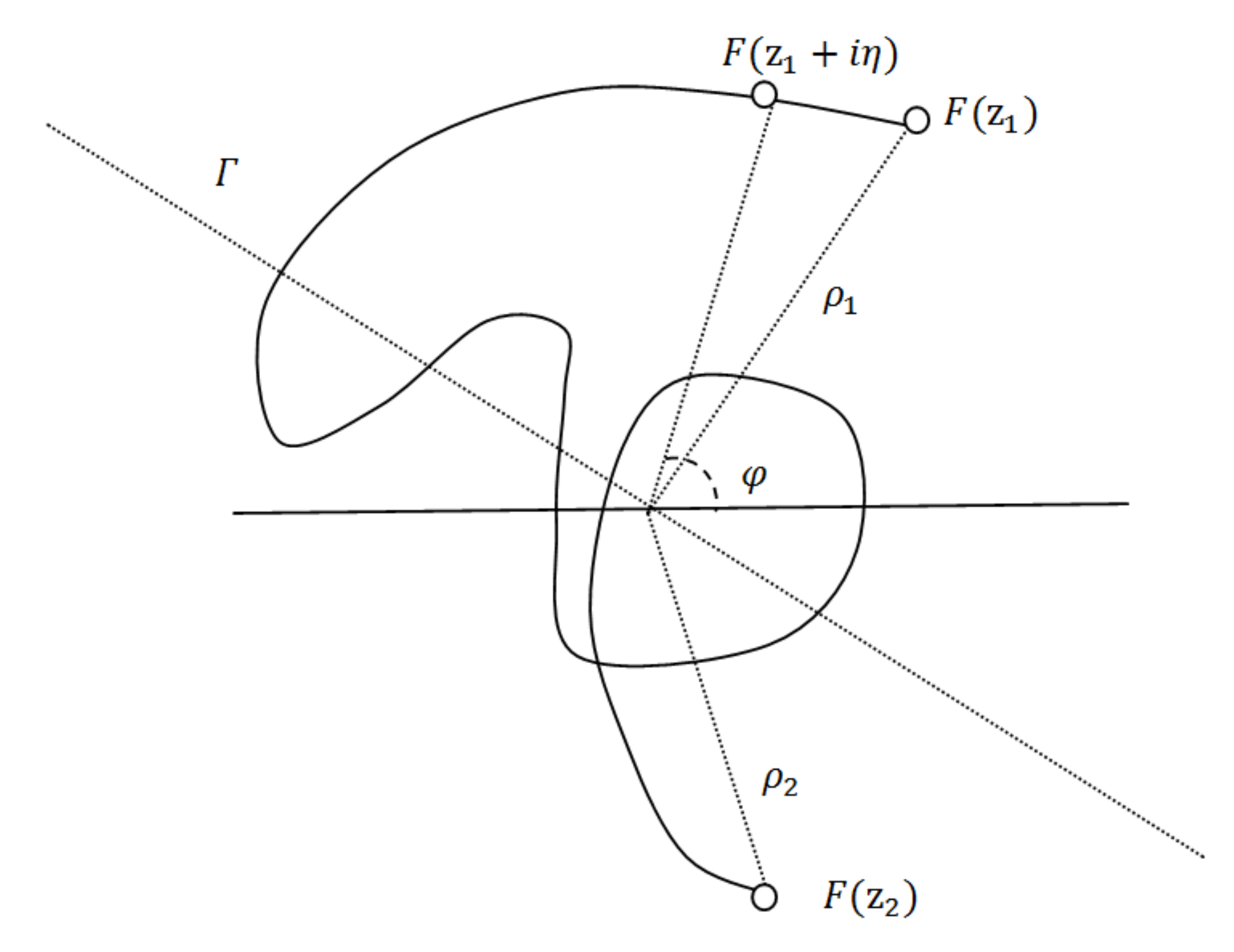}
    \end{center}
    \caption{The curve $C$ from $F(z_1)$ to $F(z_2)$.}\label{FBC.pdf}
    \end{figure}

Writing $F(z_1)=\rho_1e^{i\varphi_{z_1}}$ and $F(z_1+i\eta)=\rho_\eta e^{i\varphi_\eta}$,
where $0\leq\eta\leq T$, we have
	$$
	\Phi(\eta)=\frac{\rho_\eta}{2\rho_1}\left(e^{i(\varphi_\eta-\varphi_{z_1})}
    +e^{i(\varphi_{z_1}-\varphi_\eta)}\right)=\frac{\rho_\eta}{\rho_1}\cos(\varphi_\eta-\varphi_{z_1}).
	$$
The function $F(z)$ has no zeros on the interval $[z_1,z_2]$ by assumption. Therefore the curve
$C$ in Figure~\ref{FBC.pdf} does not pass through the origin. Consequently,
$\rho_\eta\neq 0$, so that the zeros of $\Phi(\eta)$ are completely determined by the cosine function. Thus
	\begin{equation*}
	\varphi_\eta-\varphi_{z_1}=(2n+1)\frac{\pi}{2}
	\end{equation*}
for some $n\in\Z$. Figure \ref{FBC.pdf} shows that this formula holds only at the points where
the curve $C$ intersects with the line $\Gamma$. Further, the number of zeros of $\Phi(\eta)$ on the
interval $[0,T]$, say $N$, is equal to the number of these intersection points. Then clearly
    \begin{equation}\label{less}
    N\leq n(T,\Phi(z))=:M.
    \end{equation}
Let $\{a_n\}$ denote the zero sequence of $\Phi(z)$ organized according to increasing modulus.
Since $\Phi(0)=1$, we have $a_n\neq 0$. Fix $R>T$. Then
    \begin{eqnarray*}
    n(T,\Phi(z))\log \frac{R}{T} &=& \log\frac{R^M}{T^M}\leq \log \frac{R^M}{|a_1\cdots a_M|}
    =\sum_{|a_n|\leq T}\log \frac{R}{|a_n|}\\
    &\leq & \sum_{|a_n|<R}\log\frac{R}{|a_n|}=\frac{1}{2\pi}\int_0^{2\pi}\log |F(Re^{i\theta})|\, d\theta,
    \end{eqnarray*}
where the last identity is Jensen's equality combined with $\Phi(0)=1$. This gives us
	$$
	n(T,\Phi(z))\leq \frac{\log M(R,\Phi(z))}{\log(R/T)}.
	$$
From \eqref{Phi.eq} we deduce that
	\begin{eqnarray*}
	M(R,\Phi(z)) &\leq&
    \frac{1}{2}\left(\max_{|\zeta|\leq R}\left|\frac{F(z_1+i\zeta)}{F(z_1)}\right|
    +\max_{|\zeta|\leq R}\left|\frac{\overline{F}(\overline{z}_1-i\zeta)}{\overline{F}(\overline{z}_1)}\right|\right)\\
    &=&\frac{1}{2|F(z_1)|}\left(\max_{|\zeta|\leq R}\left|F(z_1+i\zeta)\right|
    +\max_{|\zeta|\leq R}\left|F(z_1+i\overline{\zeta})\right|\right).
	\end{eqnarray*}
We have from the previous two inequalities that
	\begin{equation}\label{nR1.eq}
	n(T, \Phi(z))\leq \frac{\max_{|\zeta|\leq R}\log|F(z_1+i\zeta)|-\log |F(z_1)|}{\log (R/T)}.
	\end{equation}
The assertion follows by combining \eqref{Im.eq}, \eqref{number.eq}, \eqref{less} and
\eqref{nR1.eq}.
\end{proof}


\section{Proof of Theorem~\ref{Main.theorem}}\label{proof-sec}

The first assertion is a simple analogue of \eqref{zerofree.eq}. In order to prove \eqref{n.eq}, we use the argument principle, according to which the number of zeros of $f$ in any simply connected domain $D$ bounded by a piecewise smooth positively oriented curve $\Gamma$ is obtained from the integral
    $$
    \frac{1}{2\pi i}\int_\Gamma\frac{f'(z)}{f(z)}\,dz,
    $$
provided that $f$ has no zeros on $\Gamma$. In our case $\Gamma$ bounds a
rectangle with sides parallel to the coordinate axes. On the horizontal sides of $\Gamma$ we will use Backlund's lemma to estimate the logarithmic derivative $f'(z)/f(z)$, while on the vertical sides of $\Gamma$ we will find an asymptotic representation for $f'(z)/f(z)$.

Let $\textnormal{Re}\,(z)=x_1$ and $\textnormal{Re}\,(z)=x_2$ be two vertical lines lying in the middle of the zero-free regions of $f$ determined
by the dominant terms $H_je^{w_jz}$ and $H_ke^{w_kz}$, respectively.
Since the zero set $\{z_n\}$ of $f$ is countable, the vertical lines
$\textnormal{Im}\,(z)=y_1$ and $\textnormal{Im}\,(z)=y_2=y_1+r$ contain
no points $z_n$, except possibly for countably many $y_1\in\R$ and
for countably many $r>0$. From now on we suppose that $y_1$ and $r$ are
appropriately chosen. The idea is to allow $r$ to be arbitrarily large.

\begin{figure}[H]\label{path}
     \begin{center}
     \begin{tikzpicture}
    \draw[->](2,0)--(11,0)node[left,below,font=\tiny]{$x$};
    \draw[-,dashed](5,-1)--(5,4.5);
    \draw[-,dashed](7,-1)--(7,4.5);
    \draw[thin,blue](6,2) circle [radius=1pt];
    \draw[thin,blue](6.2,2.1) circle [radius=1pt];
        \draw[thin,blue](6.5,2.5) circle [radius=1pt];
    \draw[thin,blue](5.2,2.1) circle [radius=1pt];
        \draw[thin,blue](5.7,2.7) circle [radius=1pt];
    \draw[thin,blue](6.8,1.9) circle [radius=1pt];
        \draw[thin,blue](6.7,2.6) circle [radius=1pt];
    \draw[thin,blue](5.2,2.9) circle [radius=1pt];
        \draw[thin,blue](5.6,1.7) circle [radius=1pt];
    \draw[thin,blue](6.4,2.4) circle [radius=1pt];
       \draw[thin,blue](6,1.2) circle [radius=1pt];
    \draw[thin,blue](6.2,1.5) circle [radius=1pt];
        \draw[thin,blue](6.5,1.4) circle [radius=1pt];
    \draw[thin,blue](5.2,1.7) circle [radius=1pt];
        \draw[thin,blue](5.7,2.2) circle [radius=1pt];
    \draw[thin,blue](6.3,2.3) circle [radius=1pt];
        \draw[thin,blue](6.7,2.2) circle [radius=1pt];
    \draw[thin,blue](5.5,1.3) circle [radius=1pt];
        \draw[thin,blue](5.6,2.2) circle [radius=1pt];
    \draw[thin,blue](6.9,1.4) circle [radius=1pt];
    \draw[->](8,1)--(8,2) node[right,font=\tiny]{$\Gamma_2:z=x_2+iy$};
    \draw[-](8,2)--(8,3);
    \draw[->](4,1)--(6,1) node[below,font=\tiny]{$\Gamma_1:z=x+iy_1$} ;
    \draw[-](6,1)--(8,1);
    \draw[->](8,3)--(6,3) node[above,font=\tiny]{$\Gamma_3:z=x+iy_2$};
     \draw[-](6,3)--(4,3);
    \draw[->](4,3)--(4,2) node[left,font=\tiny]{$\Gamma_4:z=x_1+iy$};
    \draw[-](4,2)--(4,1);
    \end{tikzpicture}
     \end{center}
      \caption{The curve $\Gamma$.}\label{gamma}
    \end{figure}
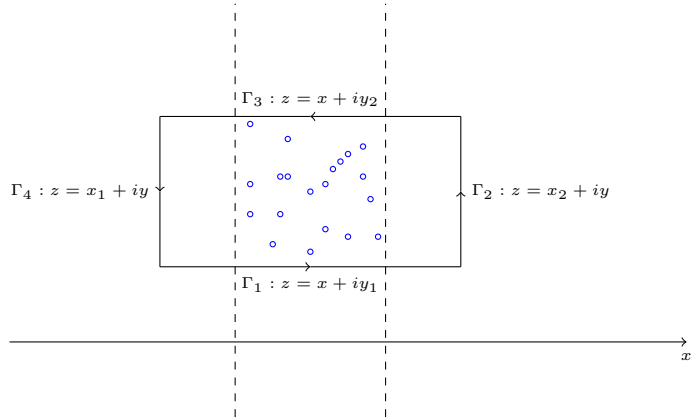

Now, we choose $\Gamma$ to be the closed Jordan curve $\Gamma:=\Gamma_1+\Gamma_2+\Gamma_3+\Gamma_4$, where
    \begin{equation*}
    \begin{array}{rll}
    \Gamma_1 &: z=x+iy_1,\ &\text{where}~x~\text{goes from}~x_1~\text{to}~x_2,\\
    \Gamma_2 &: z=x_2+iy,\ &\text{where}~y~\text{goes from}~y_1~\text{to}~y_2,\\
    \Gamma_3 &: z=x+iy_2,\ &\text{where}~x~\text{goes from}~x_2~\text{to}~x_1,\\
    \Gamma_4 &: z=x_1+iy,\ &\text{where}~y~\text{goes from}~y_2~\text{to}~y_1.
    \end{array}
    \end{equation*}
It is clear that $f$ has no zeros along $\Gamma$ by our construction, see Figure~\ref{gamma}.

We find from Backlund's lemma that
    $$
    \textnormal{Re} \left\{\frac{1}{2\pi i}\int_{\Gamma_3}\frac{f'(z)}{f(z)}\,dz\right\}
    \leq \frac{\max_{|\zeta|\leq R'}
    \log|f(x_2+iy_2+\zeta)|-\log|f(x_2+iy_2)|}
    {2(\log R'-\log|x_2-x_1|)}+\frac{1}{2},
    $$
where $R'>|x_2-x_1|$. Since
    \begin{equation*}
    \begin{split}
    |f(x_2+iy_2)|&=|1+H_1e^{w_1(x_2+iy_2)}+\cdots+H_ne^{w_n(x_2+iy_2)}|\\
    &\leq2|H_ke^{w_k(x_2+iy_2)}|=2|H_k|e^{w_kx_2}\\
    |f(x_2+iy_2)|&\geq \left|H_ke^{w_k(x_2+iy_2)}\right|-
    \sum_{l\neq k}\left|H_le^{w_l(x_2+iy_2)}\right|\\
    &= |H_k|e^{w_kx_2}-\sum_{l\neq k}|H_l|e^{w_lx_2}>0
    \end{split}
    \end{equation*}
and
	\begin{equation*}
	|f(x_2+iy_2+\zeta)| \leq 1+\sum_{l=1}^n \left|H_le^{w_l(x_2+iy_2+\zeta)}\right|
	\leq 1+\sum_{l=1}^n |H_l|e^{w_l(x_2+R')},
	\end{equation*}
we may choose $R'=2|x_2-x_1|$ to obtain
    $$
    \textnormal{Re} \left\{\frac{1}{2\pi i}\int_{\Gamma_3}\frac{f'(z)}{f(z)}\,dz\right\}=O(1).
    $$
Similarly, the corresponding integral over $\Gamma_1$ is also bounded.
Moreover,
    \begin{equation*}
    \begin{split}
    \frac{1}{2\pi i}\int_{\Gamma_2}\frac{f'(z)}{f(z)}\,dz
    &=\frac{1}{2\pi i}\left(\log f(x_2+iy_2)-\log f(x_2+iy_1)\right)\\
    &=\frac{1}{2\pi i}\bigg(\log\left(H_ke^{w_k(x_2+iy_2)}\right)\\
    &\qquad +\log\frac{1+H_1e^{w_1(x_2+iy_2)}+\cdots+H_ne^{w_n(x_2+iy_2)}}{H_ke^{w_k(x_2+iy_2)}}\bigg)+O(1)\\
    &=\frac{1}{2\pi i}\log\left(H_ke^{w_k(x_2+iy_2)}\right)+O(1)
    = \frac{w_k}{2\pi}r+O(1).
    \end{split}
    \end{equation*}
Similarly, integrating along $\Gamma_4$ but keeping in mind that the direction of integration is now reversed, it follows that
    $$
    \frac{1}{2\pi i}\int_{\Gamma_4}\frac{f'(z)}{f(z)}\,dz
    =-\frac{w_j}{2\pi}r+O(1).
    $$
From the discussion above,
    \begin{equation*}
    \begin{split}
    n(r,\Lambda_{jk})&=\frac{1}{2\pi i}\int_{\Gamma}\frac{f'(z)}{f(z)}\,dz
    =\sum_{l=1}^4 \textnormal{Re}\,\left\{\frac{1}{2\pi i}\int_{\Gamma_l}\frac{f'(z)}{f(z)}\,dz\right\}\\
    &=\frac{|w_j-w_k|}{2\pi}r+O(1),
    \end{split}
    \end{equation*}
where we have used the fact that the counting functions are non-negative.
We have proved this formula for all $r>0$ with at most countably many possible exceptions.
The assertion follows from piecewise continuity of counting functions.
Indeed, all counting functions are step functions.

\bigskip

\medskip
\noindent
\emph{J.~M.~Heittokangas}\\
\textsc{University of Eastern Finland, Department of Physics and Mathematics,
P.O.~Box 111, 80101 Joensuu, Finland}\\
\texttt{email:janne.heittokangas@uef.fi}

\medskip
\noindent
\emph{Z.-T.~Wen}\\
\textsc{Shantou University, Department of Mathematics, Daxue Road No. 243, Shantou 515063, China}\\
\texttt{e-mail:zhtwen@stu.edu.cn}


\begin{thebibliography}{99}

\bibitem{Backlund}
	Backlund R.~J.,
 	\emph{Sur les z\'eros de la fonction $\zeta(s)$ de Riemann}.
  	Comptes Rendus de I'Acad\'emie des Sciences
  	\textbf{158} (1914), 1979--1981.

\bibitem{HITW}
	Heittokangas J., K.~Ishizaki, K.~Tohge and Z.-T.~Wen,
	\emph{Zero distribution and division results for exponential polynomials}.
	Israel J.~Math.~\textbf{227} (2018), 397--421.

\bibitem{Heittokangas-Wen}
  Heittokangas J.~and ~Wen Z.-T.,
  \emph{Generalization of P\'{o}lya's zero distribution theory for exponential polynomials, plus
  sharp results for asymptotic growth}. Submitted.
  Available at arXiv: http://arxiv.org/abs/1905.08919

  \bibitem{Mora-Sepulcre-Vidal}
  Mora G., J.~M.~Sepulcre, and T.~Vidal,
  \emph{On the existence of exponential polynomials with prefixed gaps}.
  Bull.~Lond.~Math.~Soc.~\textbf{45} (2013), no.~6, 1148--1162.

 \bibitem{Moreno}
 Moreno C.~J.,
 \emph{The zeros of exponential polynomials}. I.
 Compositio Math.~\textbf{26} (1973), 69--78.

    \bibitem{Langer}
	Langer R.~E.,
	\emph{On the zeros of exponential sums and integrals}.
	Bull.~Amer.~Math.~Soc.~\textbf{37} (1931), no.~4, 213--239.

\bibitem{Poly}
    P\'olya G.,
    \emph{Geometrisches \"uber die Verteilung der Nullstellen
    spezieller ganzer Funktionen}.
    Sitz.-Ber.~Bayer.~Akad.~Wiss.~(1920), 285--290.

\bibitem{Poly2}
    P\'olya G.,
    \emph{Untersuchungen \"uber L\"ucken und Singularit\"aten
    von Potenzreihen}.
    Math.~Z.~\textbf{29} (1929), no.~1, 549--640.

\bibitem{Sch}
    Schwengeler E.,
    \emph{Geometrisches \"uber die Verteilung der Nullstellen spezieller
    ganzer Funktionen (Exponentialsummen)}.
    Diss.~Z\"urich, 1925.

\bibitem{Tamarkin}
Tamarkin J.,
\emph{Some general problems of the theory of ordinary linear differential equations and expansion of an arbitrary function in series of fundamental functions}.
Math.~Z.~\textbf{27} (1928), no.~1, 1--54.

\bibitem{Wilder}
Wilder C.~E.,
\emph{Expansion problems of ordinary linear differential equations with auxiliary conditions at more than two points.}
Trans.~Amer.~Math.~Soc.~\textbf{18} (1917), no.~4, 415--442.

\end{thebibliography}
\end{document}